\documentclass[aap]{imsart}
\usepackage{graphicx}
\RequirePackage{amsthm,amsmath,amsfonts,amssymb}
\RequirePackage[numbers]{natbib}
\usepackage{tabularx}
\usepackage{caption}
\usepackage{subcaption}
\usepackage{xspace}

\renewcommand{\tablename}{T\small{ABLE}}

\theoremstyle{plain}
\newtheorem{thm}{Theorem}
\newcommand{\MI}{\ensuremath{MI}\xspace}
\newcommand{\pvalue}{\ensuremath{p\text{-value}}}

\endlocaldefs
\begin{document}
\setlength{\abovedisplayskip}{5pt} % 上部のマージン
\setlength{\belowdisplayskip}{5pt} % 下部のマージン

\begin{frontmatter}
\title{Logarithmic Asymptotic Relations Between $p$-Values and Mutual Information}
\runtitle{Asymptotic $p$-Value--MI Relations}

\begin{aug}
%%%%%%%%%%%%%%%%%%%%%%%%%%%%%%%%%%%%%%%%%%%%%%%
%% Only one address is permitted per author. %%
%% Only division, organization and e-mail is %%
%% included in the address.                  %%
%% Additional information can be included in %%
\author[A]{\fnms{Tsutomu} \snm{Mori}\ead[label=e1]{tmori@fmu.ac.jp}}
\and
\author[B]{\fnms{Takashi} \snm{Kawamura}\ead[label=e3]{tkawamur@fmu.ac.jp}}
%%%%%%%%%%%%%%%%%%%%%%%%%%%%%%%%%%%%%%%%%%%%%%
%% Addresses                                %%
%%%%%%%%%%%%%%%%%%%%%%%%%%%%%%%%%%%%%%%%%%%%%%
\address[A]{Department of Human Lifesciences, Fukushima Medical University School of Nursing \printead{e1}}

\address[B]{Department of Human Lifesciences, Fukushima Medical University School of Nursing \printead{e3}}
\end{aug}

\begin{abstract}
We establish a precise connection between statistical significance in dependence testing and information-theoretic dependence as quantified by Shannon mutual information (MI).
In the absence of prior distributional information, we consider a maximum-entropy model and show that the probability associated with the realization of a given magnitude of MI takes an exponential form, yielding a corresponding tail-probability interpretation of a $p$-value.
In contingency tables with fixed marginal frequencies, we analyze Fisher's exact test and prove that its $p$-value $P_F$ satisfies a logarithmic asymptotic relation of the form $MI=-(1/N)\log P_F + O(\log(N+1)/N)$ as the sample size $N\to\infty$.
These results clarify the role of MI as the exponential rate governing the asymptotic behavior of $p$-values in the settings studied here, and they enable principled comparisons of dependence across datasets with different sample sizes.
We further discuss implications for combining evidence across studies via meta-analysis, allowing mutual information and its statistical significance to be integrated in a unified framework.
\end{abstract}



\end{frontmatter}

%%%%%%%%%%%%%%%%%%%%%%%%%%%%%%%%%%%%%%%%%%%%%%
%% Please use \tableofcontents for articles %%
%% with 50 pages and more                   %%
%%%%%%%%%%%%%%%%%%%%%%%%%%%%%%%%%%%%%%%%%%%%%%
%\tableofcontents

%%%%%%%%%%%%%%%%%%%%%%%%%%%%%%%%%%%%%%%%%%%%%%
%%%% Main text entry area:
\section{Introduction}
Despite their close conceptual relationship, probability theory and information theory have long developed as largely independent frameworks.
For a single random variable, Shannon's self-information $I=-\log p$ \cite{r16} provides an exact and well-established correspondence between probability and information.
For two random variables, however, no universally accepted information-theoretic quantity has been identified that directly corresponds to the statistical significance probability, commonly expressed as a $p$-value.
Clarifying this relationship is not only of conceptual interest but also of practical importance in modern data analysis, where large-scale and high-dimensional dependence structures are routinely investigated.

In this paper, we revisit the relationship between probability theory and statistics and information theory through the lens of mutual information (MI).
Rather than asserting a literal identity between MI and the $p$-value, we formulate and prove precise statements that clarify \emph{when} and \emph{in what sense} a $p$-value admits an information-theoretic representation.
The central theme is that, in the settings considered here and under appropriate asymptotic conditions, the logarithm of a $p$-value admits an information-theoretic characterization through an exponential-rate form.
This perspective highlights an exponential-rate structure in dependence testing and clarifies how statistical significance relates to information-theoretic dependence.
Throughout, the term ``equivalence'' is used in this logarithmic asymptotic sense, rather than as an exact finite-sample identity.

We consider two complementary settings.
First, when no prior information about the distributions of the random variables is available, the principle of maximum entropy \cite{r9} provides a canonical probabilistic model.
In this case, we show that the probability associated with the realization of a given magnitude of mutual information follows an exponential form, and the corresponding tail probability yields a $p$-value that is exponentially determined by MI.
Second, when prior information is available in the form of fixed marginal frequencies, as in contingency tables, we analyze Fisher's exact test \cite{r6}.
In this classical statistical setting, we prove that the Fisher $p$-value is asymptotically related to MI through a precise scaling with the sample size.
More precisely, Theorems~1 and~2 establish an exponential form for the probability associated with mutual information under a maximum-entropy model, yielding a corresponding tail-probability interpretation of the $p$-value.
Theorem~3 then shows that, for contingency tables with fixed margins, the Fisher exact $p$-value ($P_F$) satisfies an asymptotic relation of the form $MI=-(1/N)\log P_F + O\!\left(\frac{\log(N+1)}{N}\right)$ as $N\to\infty$, where log is natural.
We emphasize that Theorem~3 is the main theorem of the paper: it concerns the exact Fisher $p$-value, defined as a tail sum under fixed margins. By contrast, Theorems~1--2 provide a supporting, model-based calibration under a maximum-entropy baseline.
Mutual information is closely connected to likelihood-ratio statistics and to large-deviation principles: in multinomial models it equals the Kullback-Leibler divergence between the empirical joint distribution and the product of its empirical marginals, and the log-likelihood ratio for independence is $N\cdot MI$ (equivalently, the $G^2$ statistic is $2N\cdot MI$) \cite{r4,r25}.
Sanov-type large deviations interpret this divergence as an exponential rate and motivate chi-square or large-deviation approximations to $p$-values (see, e.g., \cite{r24,r23}). Our contribution is complementary: Theorem~3 treats the exact two-sided Fisher $p$-value under fixed margins and proves that its exponential rate is governed by $MI$, with explicit finite-$N$ bounds on the logarithmic scale that quantify the tail-sum contribution.
We use $MI$ to denote the mutual information throughout the paper.
Although these two settings are conceptually distinct, their asymptotic consequences coincide, revealing a unified structure underlying probabilistic and information-theoretic measures of dependence.

Meanwhile, in information theory, Shannon's MI, which represents the information exchanged between two random variables, is excellent among the measures of dependence between variables \cite{r4, r11, r20}. MI is unique in its close ties to Shannon entropy \cite{r11}, that is, $MI=H(X)+H(Y)-H(X,Y)$, using entropy $H(X)$ and $H(Y)$ of random variables $X$ and $Y$, respectively, and their joint entropy $H(X,Y)$. MI has advantages over the correlation coefficient because it measures types of dependence other than linear dependence \cite{r2}. Another advantage of MI is that it is zero if and only if the two random variables are strictly independent \cite{r11}. These properties demonstrate that MI is an orthodox measure in information theory \cite{r10}.

In the case of two discrete random variables, their MI can be defined as follows. Assume that $X$ and $Y$ take values from $X_1$ to $X_m$ and from $Y_1$ to $Y_n$, respectively. Let $p(X_i, Y_j)$ be the joint probability that $X$ takes $X_i$ and $Y$ takes $Y_j$. Let $p(X_i)$ be the marginal probability that $X$ takes $X_i$ and $p(Y_j)$ be the marginal probability that $Y$ takes $Y_j$. Then, the MI of $X$ and $Y$, $MI\ge 0$, is defined as
\begin{equation}
MI= \sum_{i=1}^{m}\sum_{j=1}^{n} p(X_i, Y_j)\log{\frac{p(X_i, Y_j)}{p(X_i)p(Y_j)}}.
\end{equation}

\subsection{Posing a problem}
 We aim to solve the following problems related to the mathematical foundation that bridges probability statistics and information theory. Considering the completeness of both mathematical systems, we might admit that there is little room for advancement. However, a knowledge gap exists in the theoretical area that should connect both systems. We describe the gap as the following two concerns regarding the measurement of interdependence between two random variables in terms of each mathematical system.

\begin{enumerate}
  \item[a)]The $p$-value, a measure in probability statistics, is not quantitative because it varies largely depending on the sample size.

  \item[b)]MI, a measure in information theory, has been studied in depth, whereas its probability-statistical characteristics remain unclear.
\end{enumerate}

 Interdependence between random variables is a central theme in both theoretical systems. Since both the $p$-value in hypothesis testing and mutual information aim to quantify departures from independence, one expects a connection between them. While links between likelihood-based statistics and information measures are well known, the connection between mutual information and \emph{exact} significance measures---such as Fisher's exact-test $p$-value under fixed margins---appears comparatively less explicit, especially with finite-$N$ quantitative control on the logarithmic scale. Motivated by this gap, in 1.2 for probability theory and 1.3 for statistics we state our Theorems \ref{th1}-\ref{th3}.
 
\subsection{Maximum-entropy calibration in probability theory}
First, by applying probability theory, we calculate the probability that information exchange arises between two random variables with unknown distributions. 
 
\subsubsection{Principle of maximum entropy}
Let ${X}$ and ${Y}$ be two discrete random variables that take a finite number of nonnegative rational values with no prior information about distributions. We apply the principle of maximum entropy \cite{r9} to $X$ and $Y$. According to this, we assume that both $X$ and $Y$ follow the uniform distribution.

\subsubsection{Principle of equal probability in statistics}
When two random variables $X$ and $Y$ follow the uniform distribution, each value that $X$ and $Y$ can take occurs with equal probability. This property corresponds to the principle of equal probability in statistical mechanics, which demonstrates the close relationship between physics and probability statistics. Indeed, statistical mechanics has been regarded as a form of statistical inference rather than a physical theory \cite{r9}. Therefore, in terms of statistical mechanics, let $W_X$ and $W_Y$ be the number of states that $X$ and $Y$ can take, respectively. Then their information entropies $H(X)$ and $H(Y)$ are represented as
\begin{equation}
H(X)=\log{W_X}, \qquad H(Y)=\log{W_Y},
\end{equation}
respectively. These correspond to Boltzmann's principle, which represents thermodynamic entropy as $S=k_B\log{W}$, where $k_B$ is the Boltzmann constant and $W$ is the number of microscopic states.

Moreover, we also apply the maximal entropy principle to the joint random variable $(X, Y)$ and assume its uniform distribution. Then the joint information entropy $H(X, Y)$ is represented as
\begin{equation}
H(X, Y)= \log{W_{XY}},
\end{equation}
where $W_{XY}$ is the number of states of the joint random variable $(X, Y)$.

\subsubsection{An exponential-form relation between information and probability}
With respect to the probability of MI, we obtain the following Theorem \ref{th1}.
\vspace{10pt}

\begin{thm}\label {th1}
Let $X$ and $Y$ be random variables that follow a uniform distribution. Then the probability ${P_{MI}}$ that the magnitude of MI shared by them becomes $MI$ is represented as
\begin{equation}
P_{MI} = e^{-MI}, \qquad MI = -\log{P_{MI}}.
\end{equation}
\end{thm}

\vspace{10pt}

These formulas demonstrate that $MI$ is an alternative expression of $P_{MI}$, the probability of exchanging information. These formulas provide an exponential-form relation between mutual information and its realization probability under the maximum-entropy baseline. Moreover, (4) is isomorphic to the formulas for the realization probability $P_I$ of self-information $I$,
\begin{equation}
P_I = e^{-I}, \qquad I = -\log{P_I}.
\end{equation}
Therefore, the fact that information and its realization probability are represented by the logarithmic and exponential functions of each other holds not only for self-information $I$ of a single random variable but also for MI shared by two random variables.

In association with $I$ and $MI$, the thermodynamic entropy $S$ in the canonical ensemble of statistical mechanics is known to follow the exponential distribution. This is an important principle of physics and is called canonical distribution. Considering the generality of the relationship between information and its realization probability represented by (4) and (5), we term these probability distributions the infocanonical distributions \cite{r13}.

\subsubsection{The $p$-value based on information theory} In (4), $P_{MI}$ represents the probability mass function that MI shared by two discrete random variables $X$ and $Y$ is $MI$. However, when the sample spaces of $X$ and $Y$ are sufficiently large, we obtain the probability density function 
\begin{equation}
f_{MI} = e^{-MI},
\end{equation}
by continuing and normalizing $P_{MI} = e^{-MI}$. Then, by taking the limit $MI \rightarrow \infty$ under the condition that the entropies $H(X)$ and $H(Y)$ tend to infinity, we calculate the $p$-value between $X$ and $Y$ as
\begin{equation}
\pvalue = \int_{t=MI}^{\infty} e^{-t} d t = e^{-MI}.
\end{equation}
Therefore, $P_{MI} = e^{-MI}$ represents not only the probability that the magnitude of MI is $MI$ but also the $p$-value, which is the significance probability that the magnitude of MI is greater than or equal to $MI$. If $p$-value = 0.05, then $MI = 2.9957 \simeq 3.00$. Hence, if $MI \ge 3.00$, then information exchange is not accidental, and $X$ and $Y$ have significant interdependence.

As mentioned above, this exponential-form relation in the maximum-entropy setting connects $MI$, its realization probability $P_{MI}$, and its $p$-value. Thus, in the field of probability theory, we have succeeded in unifying the $p$-value and \MI, which are distinct measures in probability theory and information theory, respectively.

\subsection{Logarithmic asymptotics in statistics: Fisher's exact test}
Next, we consider how a corresponding logarithmic relation arises in statistics,

\subsubsection{Occurring probability of information exchange per one observation}
To estimate the interdependence  between two random variables $X$ and $Y$, we consider a trial in which information exchange is observed between them repeatedly. Specifically, this corresponds to observing the states of $X$ and $Y$ in a repeated manner and summing up information at the end. When information exchange arises, if we assume that the MI of magnitude $MI$ arises per one observation, then the probability that information exchange arises is $P_{MI} = e^{-MI}$ for one observation.

\subsubsection{Occurring probability of information exchange during many observations}
When we repeat the above trial $N$ times, we obtain the following Theorem \ref{th2}.
\vspace{10pt}

\begin{thm}\label{th2}
Let $X$ and $Y$ be random variables that follow a uniform distribution. If we repeat the trial of observing information exchange between them $N$ times, then the probability $P_{N\cdot MI}$ that the total magnitude of \MI shared by them becomes $N\cdot MI$ is represented as
\begin{equation}
P_{N\cdot MI} = e^{-N\cdot MI}, \qquad MI = -\frac{1}{N}\log{P_{N\cdot MI}}.
\end{equation}
\end{thm}
\vspace{10pt}

\subsubsection{The $p$-value based on information theory}
In (8), $P_{N\cdot MI}$ represents the probability mass function such that the sum of MI shared by two discrete random variables $X$ and $Y$ is $N\cdot MI$. However, when the sample spaces of $X$ and $Y$ are sufficiently large, we obtain the probability density function
\begin{equation}
f_{N\cdot MI} = Ne^{-N\cdot MI}
\end{equation}
by continuing and normalizing $P_{N\cdot MI} = e^{-N\cdot MI}$. Then, by taking the limit of the integral variable $MI \rightarrow \infty$ under the condition that the entropies $H(X)$ and $H(Y)$ go to infinity, we calculate the $p$-value between $X$ and $Y$ as
\begin{equation}
\pvalue = \int_{t=MI}^{\infty} Ne^{-Nt} dt = e^{-N\cdot MI}.
\end{equation}

\setlength\textfloatsep{22pt}

\begin{figure}
\centering
\includegraphics[width=7cm]{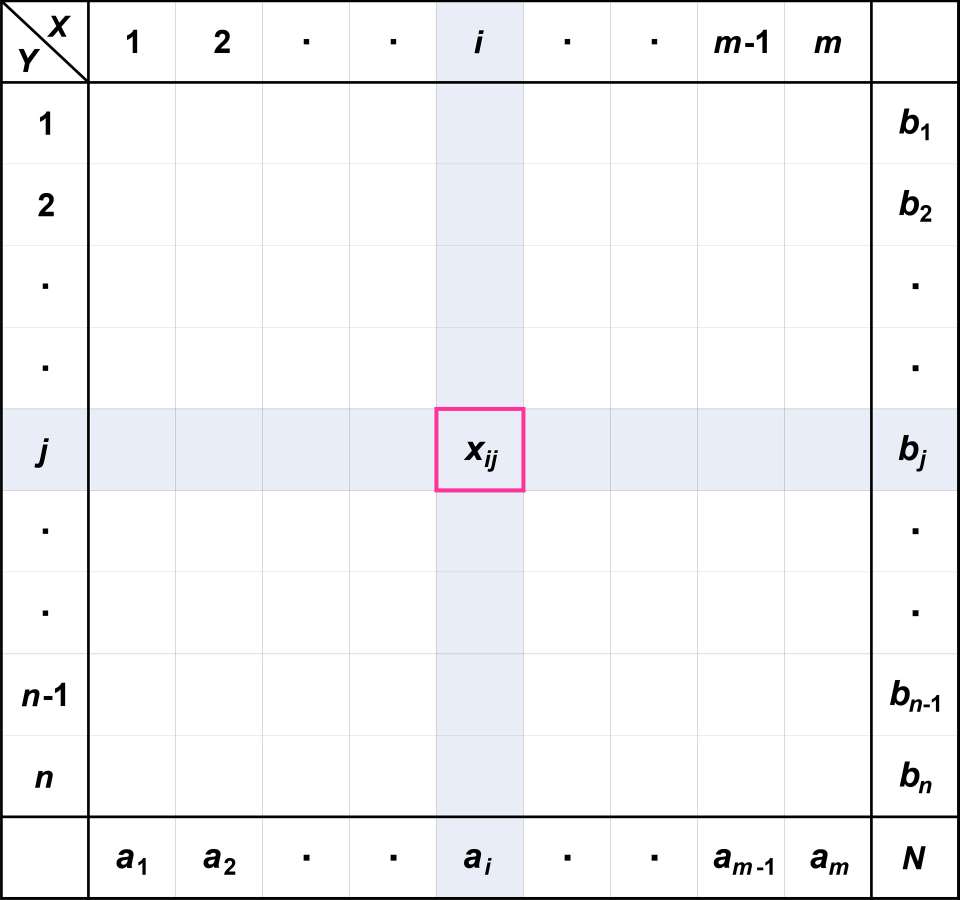}
\caption{\textit{$m\times n$ contingency table when random variables $X$ and $Y$ take $X_1$-$X_m$ and $Y_1$-$Y_n$, respectively. $x_{ij}$ is the joint frequency, $a_1$-$a_m$ and $b_1$-$b_n$ are marginal frequencies, and $N$ is the sample size.}}
\label{Figure 1}
\end{figure}

Thus, $P_{N\cdot MI} = e^{-N\cdot MI}$ represents not only the probability that the total MI of $N\cdot MI$ is realized but also the $p$-value, which is the significance probability that MI greater than or equal to $MI$ arises. If $p$-value = 0.05, then $MI = 2.9957 / N \simeq 3.00 / N$. Hence, if $MI \ge 3.00 / N$, then information exchange is not accidental, and $X$ and $Y$ have significant interdependence.

\subsubsection{Analysis using conventional statistics}
Next, we assign the results of the above trials to an $m\times n$ contingency table (Figure \ref{Figure 1}) \cite{r13} and analyze it according to conventional statistics. Let $x_{ij}$ be the joint frequency of each cell of the contingency table, and let $a_i$ and $b_j$ be the marginal frequencies. Then the probability $P_C$ that the result corresponding to the contingency table in Figure \ref{Figure 1} occurs is given by
\begin{equation}
P_C=\frac{\prod_{i=1}^m a_i! \prod_{j=1}^n b_j!}{N! \prod_{i, j} x_{ij}!}.
\end{equation}
The hypergeometric distribution probability $P_H$ is the probability $P_C$ that the observed result is obtained in the contingency table with known marginal frequencies. This $P_H$ is a probability mass function. Meanwhile, Fisher's exact probability $P_F$ is the sum of probabilities $P_C$ that the observed result or more unlikely results than it are obtained in the contingency table with fixed marginal frequencies and expressed as
\begin{equation}
P_F= \sum_{P_C\le P_H}P_C
\end{equation}
\vspace{1pt}

Next, we present Theorem \ref{th3}, which illustrates the information-theoretical properties of $P_F$. Because $P_F$ expresses the $p$-values of contingency tables with fixed marginal distributions, it considers a scenario opposite to case 1.2 of no prior information about the distribution of the random variables. However, even in this case 1.3, the relationship between information and probability is an important issue. We show that an analogous logarithmic asymptotic relation also holds in this setting. Theorem \ref{th3} is stated as follows.

\vspace{10pt}

\begin{thm}\label{th3}
Let $P_F$ be the $p$-value of Fisher's exact test defined as in (12) for an $m\times n$ contingency table with fixed marginal totals, where $m$ and $n\in\mathbb{N}$ are fixed and the total sample size is $N$.
Write $x_{ij}$ for the cell counts and assume that, along a sequence of such tables with $N\to\infty$, the proportions $x_{ij}/N$ converge to a limit $p_{ij}$ with $p_{ij}>0$ for all $i,j$.
Let $MI$ denote the mutual information computed from the empirical proportions $x_{ij}/N$.
Then, as $N\to\infty$,
\begin{equation}
-\frac{1}{N}\log P_F = MI + O\!\left(\frac{\log(N+1)}{N}\right).
\end{equation}
\end{thm}
\vspace{10pt}

The assumption is satisfied, for example, when the contingency tables arise as empirical counts from repeated sampling from a fixed joint distribution with full support.
We view Theorem \ref{th3} as providing a logarithmic asymptotic link between $P_F$ and MI. This is because this formula demonstrates its fundamental nature by connecting the two authentic measures of interdependence in probability theory and information theory. Using our formula, we can calculate $MI$ from $P_F$, which allows us to use the merit of not only $P_F$ but also $MI$. For example, $P_F$ has the drawback of being sensitive to the sample size \cite{r14}, whereas $MI$ converted from $P_F$ can specify interdependence irrespective of the sample size.
In addition to Theorem \ref{th3}, when sample size $N$ goes to infinity, the following formulas hold concerning $MI$, $P_H$ and $P_F$:
\begin{equation}
-\frac{1}{N}\log P_H = MI + O\!\left(\frac{\log(N+1)}{N}\right),\qquad -\frac{1}{N}\log P_F = MI + O\!\left(\frac{\log(N+1)}{N}\right),
\end{equation}

Moreover, in the following sections, we prove that (14) is asymptotically equal to the equation on the right of (8). As mentioned above, (8) and (10) demonstrate the logarithmic asymptotic relationship among $N\cdot MI$, the probability mass function of information exchange, $P_{N\cdot MI}$, and the $p$-value. Additionally, (14) shows the asymptotic relationship among $N\cdot MI$, the probability mass function in statistics, $P_H$, and the $p$-value, $P_F$. Together, these relations provide an information-theoretic interpretation of statistical significance in contingency-table dependence testing.
Thus, in the field of statistics, we have succeeded in unifying the $p$-value and MI, which are distinct measures of statistics and information theory, respectively. Given the generality of these measures, our theorems can be applied to various fields of science, such as medicine and biology. In particular, we exploit the theorems to precisely estimate MI between functionally interacting genes in biological systems.

The main contributions of this paper may be summarized as follows:
\begin{itemize}
\item We prove that, for contingency tables with fixed margins (with $m,n$ fixed), $P_F$ satisfies the logarithmic asymptotic relation $MI = -(1/N)\log P_F + O\!\left(\frac{\log(N+1)}{N}\right)$ as $N\to\infty$ (Theorem~\ref{th3}).
\item In the $2\times2$ case we derive explicit quantitative error bounds via a decomposition into a Stirling-approximation error and a tail-sum error, and we extend the logarithmic asymptotics to general $m\times n$ tables.
\item We discuss implications of this rate interpretation for comparing dependence strengths across datasets with different sample sizes and for combining evidence across studies via meta-analysis.
\end{itemize}

This paper is organized as follows: We prove Theorems \ref{th1} and \ref{th2} in Section 2. We prove Theorem \ref{th3} for a $2\times 2$ contingency table in Section 3 and extend it for a general $m\times n$ contingency table in Section 4. In Section 5, we explain the numerical simulations that we used to verify Theorem \ref{th3}. In Section 6, we demonstrate that Theorem \ref{th3} can be applied to the meta-analysis of MI for any dimension, which produces a low $p$-value. In Section 7, we discuss the advantages and applications of our Theorems, and illustrate a wide range of benefits that arise from making the techniques of information theory and probability statistics available.

\section{Proofs of Theorems \ref{th1} and \ref{th2}}
In this section, we prove Theorems \ref{th1} and \ref{th2}.
In this section we adopt a maximum-entropy (uniform) baseline model as a calibration device. Theorems~\ref{th1} and~\ref{th2} should be read as model-based exponential-form identities under this baseline, and not as the sampling distribution of an MI estimator. These results are logically independent of the Fisher exact-test asymptotics in Theorem~\ref{th3}.
\subsection{Proof of Theorem \ref{th1}}

Let $W_X$ and $W_Y$ be the number of states that $X$ and $Y$ can adopt under a certain condition, respectively. By contrast, the number of all states that $X$ and $Y$ can take are defined as $W_{Xall}$ and $W_{Yall}$, respectively. Then, the state probabilities $p_X$ and $p_Y$ can be expressed as $p_X = W_X / W_{Xall}$ and $p_Y = W_Y / W_{Yall}$, respectively. According to the assumption in Section 1, the random variables $X$ and $Y$ follow the uniform distribution. Therefore, the occurring probabilities $p_X$ and $p_Y$, which are proportional to their number of states, satisfy
\begin{equation}
p_X \propto W_X = \exp{[H(X)]}, \qquad p_Y \propto W_Y = \exp{[H(Y)]},
\end{equation}
respectively. Similarly, the occurring probability $p_{XY}$ of the joint random variable $(X, Y)$ satisfies
\begin{equation}
p_{XY} \propto W_{XY} = \exp{[H(X, Y)]} .
\end{equation}

If the state probability that $X$ and $Y$ are independent is one, then we calculate the probability $P_{MI}$ that the magnitude of MI shared by them becomes $MI$ as
\begin{equation}
P_{MI} = \frac{W_{XY}}{W_X W_Y} = \frac{\exp{[H(X, Y)]}}{\exp{[H(X)+H(Y)]}} = \frac{\exp{[H(X)+H(Y)-MI]}}{\exp{[H(X)+H(Y)]}} = e^{-MI}.
\end{equation}

Hence, we obtain
\begin{equation}
P_{MI} = e^{-MI}, \qquad MI = -\log{P_{MI}}, \notag
\end{equation}
which completes the proof of Theorem \ref{th1}. $\square$
\subsection{Proof of Theorem \ref{th2}}
When we repeat the above trial $N$ times independently, the expectation value of the sum of MI exchanged between $X$ and $Y$ is $N\cdot MI$. Simultaneously, its expected realization probability, $P_{N\cdot MI}$, is represented as
\begin{equation}
P_{N\cdot MI} = (e^{-MI})^N = e^{-N\cdot MI} \notag
\end{equation}
which coincides with the formula obtained by substituting $N\cdot MI$ for $MI$ in (4).
Thus, (8) holds for any $N$, which completes the proof of Theorem \ref{th2}.  $\square$

\section{Proof of Theorem \ref{th3} for a $2\times 2$ contingency table}
In this section, we prove Theorem \ref{th3} for a $2\times 2$ contingency table. We consider the following two contingency tables, Tables 1 and 2, for two random variables $A$ and $B$, which each take two values $A_1$ and $A_2$ and $B_1$ and $B_2$, respectively.

\subsection{Setup}

Table 1 has integer entries, whereas Table 2 has real number entries that include or do not include irrational numbers. These tables can be converted to each other as described below.

First, we mention the translation of Table 1 into Table 2. Table 1 shows the frequency of the combination of the variables, whereas Table 2 shows the relative frequency obtained from Table 1 by dividing by $N$. Then, $X_0$, $X_1$, $X_2$, and $X_3$ are rational numbers between 0 and 1 that approach the true joint probabilities, $p(A_1, B_1)$, $p(A_2, B_1)$, $p(A_1, B_2)$, and $p(A_2, B_2)$, respectively, as $N$ goes to infinity. 

Second, Table 2 can be transformed into the form of Table 1 by assuming a large sample size $N$. In applications, the observed data are integer counts (Table 1), and Table 2 records the corresponding empirical proportions, so $X_0,\ldots,X_3$ are rational numbers with denominator $N$. For notational convenience we occasionally view $X_0,\ldots,X_3$ as real numbers; along sequences with $N\to\infty$, such rational proportions can approximate any limiting probability vector arbitrarily closely. In the remainder of this paper, we work with a sequence of contingency tables with increasing $N$ for which $NX_0$ to $NX_3$ are integers.

Using Table 2, the MI of $A$ and $B$ is defined according to (1) as
\begin{align}
MI &= X_0\log{\frac{X_0}{(X_0+X_2)(X_0+X_1)}}+X_1\log{\frac{X_1}{(X_1+X_3)(X_0+X_1)}} \\
   &\ +X_2\log{\frac{X_2}{(X_0+X_2)(X_2+X_3)}}+X_3\log{\frac{X_3}{(X_1+X_3)(X_2+X_3)}} \notag \\
   &= \sum_{k=0}^{3} X_k\log{X_k}-(X_0+X_1)\log{(X_0+X_1)}-(X_0+X_2)\log{(X_0+X_2)}\notag \\
   &\ -(X_1+X_3)\log{(X_1+X_3)}-(X_2+X_3)\log{(X_2+X_3)},\notag 
\end{align}
where $0\log{0}$ is defined as $0$.

\newpage

\begin{center}
T{\small ABLE} 1\\
\textit{Contingency table of the observed frequency.} \\
\textit{All the values of the cells are non-negative integers.} \\
\begin{tabularx}{90mm}{cccc} \hline
\rule[-0.5em]{0em}{2em} & $A_1$ & $A_2$ & Total\\ [3pt]\hline
\rule[-0.5em]{0em}{2em} $B_1$ & $NX_0$ & $NX_1$ & $N(X_0+X_1)$ \\  [3pt]
$B_2$ & $NX_2$ & $NX_3$ & $N(X_2+X_3)$ \\ [3pt]\hline
\rule[-0.5em]{0em}{2em} Total & $N(X_0+X_2)$ & $N(X_1+X_3)$ & $N$ \\ [3pt]\hline
\end{tabularx}
\end{center}
\vspace{10pt}

\begin{center}
T{\small ABLE} 2\\
\textit{Contingency table of the relative frequency.}  \\
\textit{All the values of the cells are non-negative real numbers.}
\begin{tabularx}{90mm}{cccccccc} \hline
\rule[-0.5em]{0em}{2em} & & $A_1$ & & $A_2$ & & & Total\\ [3pt]\hline
\rule[-0.5em]{0em}{2em} $B_1$ & & $X_0$ & & $X_1$ & & & $X_0+X_1$ \\ [3pt]
$B_2$ & & $X_2$ & & $X_3$ & & &  $X_2+X_3$ \\ [3pt]\hline
\rule[-0.5em]{0em}{2em} Total & & $X_0+X_2$ & & $X_1+X_3$ & & & $1$ \\ [3pt]\hline
\end{tabularx}
\end{center}
\vspace{1pt}

\subsection{MI and hypergeometric distribution probability}

In the following subsection, we examine the relationship between $P_F$ and $MI$ for $2\times 2$ contingency tables. We divide $P_F$ into two terms as $P_F=P_H+P_{Frem}$, where $P_H$ is the main term and $P_{Frem}$ is the sum of the remaining terms. First, $P_H$ is the probability that the observed result is obtained, which is the hypergeometric distribution probability. Second, $P_{Frem}$ represents the probability that less possible results than the observed result occur. As shown below, $P_{Frem}$ asymptotically becomes negligible compared with the main term.

First, we outline the information conversion of the main term of $P_F$, that is, $P_H$. We calculate the main term as
\setlength{\abovedisplayskip}{4pt} 
\setlength{\belowdisplayskip}{14pt} 

\begin{align}
P_{H} &= \frac{\binom{N(X_0+X_1)}{NX_0}\binom{N(X_2+X_3)}{NX_2}}{\binom{N}{N(X_0+X_2)}} \\
      &= \frac{[N(X_0+X_1)!]}{[(NX_0)!(NX_1)!]}\times\frac{[N(X_2+X_3)!]}{[(NX_2)!(NX_3)!]}\div\frac{N!}{[N(X_0+X_2)]![N(X_1+X_3)]!}. \notag
\end{align}

Taking the logarithm, we obtain
\setlength{\abovedisplayskip}{7pt} 
\setlength{\belowdisplayskip}{10pt} 
\begin{align}
-\log{P_H} &= \sum_{k=0}^{3} \log{(NX_k)!}+\log{N!}-\log{[N(X_0+X_1)]!} \\
           &\ -\log{[N(X_0+X_2)]!}-\log{[N(X_1+X_3)]!}-\log{[N(X_2+X_3)]!}. \notag
\end{align}
\setlength{\abovedisplayskip}{7pt} 
\setlength{\belowdisplayskip}{8pt} 
To derive an approximate formula, we apply Stirling's formula, $\log {n!}\approx n\log {n} - n$, where $n$ is large. When $N$ is large, $NX_0$ to $NX_3$ are large, and we can apply Stirling's formula. Using $X_0 +X_1 +X_2 +X_3 = 1$ and (18), we obtain
\begin{align}
-\log{P_H} &\approx \sum_{k=0}^{3}(NX_k\log{NX_k}-NX_k)+N\log{N}-N \\
           &\ -N(X_0+X_1)\log{N(X_0+X_1)}+N(X_0+X_1) \notag \\
           &\ -N(X_0+X_2)\log{N(X_0+X_2)}+N(X_0+X_2) \notag \\
           &\ -N(X_1+X_3)\log{N(X_1+X_3)}+N(X_1+X_3) \notag \\
           &\ -N(X_2+X_3)\log{N(X_2+X_3)}+N(X_2+X_3) \notag \\
           &= \sum_{k=0}^{3}(NX_k\log{N}+NX_k\log{X_k})+N\log{N}-N(X_0+X_1)\log{N}\notag \\
           &\ -N(X_0+X_1)\log{(X_0+X_1)}-N(X_0+X_2)\log{N} \notag \\
           &\ -N(X_0+X_2)\log(X_0+X_2)-N(X_1+X_3)\log{N} \notag \\
           &\ -N(X_1+X_3)\log{(X_1+X_3)}-N(X_2+X_3)\log{N}\notag \\
           &\ -N(X_2+X_3)\log{(X_2+X_3)}\notag \\
           &= \sum_{k=0}^{3}NX_k\log{X_k}-N(X_0+X_1)\log{(X_0+X_1)}\notag \\
           &\ -N(X_0+X_2)\log{(X_0+X_2)}-N(X_1+X_3)\log{(X_1+X_3)}\notag \\
           &\ -N(X_2+X_3)\log{(X_2+X_3)}\notag \\
           &= N\cdot MI. \notag
\end{align}
\setlength{\abovedisplayskip}{7pt} 
\setlength{\belowdisplayskip}{10pt} 
Hence, if $N$ is sufficiently large and if $NX_0$ to $NX_3$ approach integers, then
\begin{equation}
MI\approx -\frac{1}{N}\log{P_H}.
\end{equation}
\vspace{3pt}
Thus, (22) indicates that $MI$ is approximately equal to the logarithm of $P_H$, divided by $N$. Given that $MI$ and $P_H$ have been defined differently, this equivalence highlights a close connection between them that is not always made explicit. Additionally, because $P_H$ represents the hypergeometric distribution, $MI$ is inherently related to sampling without replacement rather than that with replacement represented by the binomial distribution.

\vspace{5pt}
\subsection{Evaluation of the error in Theorem \ref{th3}}

Next, we evaluate the error between $MI$ and $-(\log{P_F})/N$ in (13) by multiplying both sides by $N$. The error in Theorem 3 comprises two parts. The first part is the error of Stirling's formula applied to the main term $P_H$ and the second part is the sum of the remaining terms $P_{Frem}=P_F-P_H$. To assess the errors, let $ER_1 = -\log{P_H} - N\cdot MI$ be the first part of the error and let $ER_2 = -\log{P_H} - (-\log{P_F})$ be the second part of the error. The latter is related to $P_{Frem}$.

\subsubsection{Evaluation of the first part of the error}

We estimate the first part of the error $ER_1 = -\log P_H -\ N\cdot MI$. Stirling's formula in its exact form is expressed as

\setlength{\abovedisplayskip}{12pt} 
\setlength{\belowdisplayskip}{8pt} 
\begin{align}
N! &=\Gamma(N+1) \\
           &=(N+1)^{N+\frac{1}{2}}\exp{[-(N+1)]}\sqrt{2\pi}\exp{\Bigg[\sum_{n=1}^{\infty}\frac{(-1)^{n-1}B_{2n}}{2n(2n-1)(N+1)^{2n-1}}\Bigg]}, \notag
\end{align}
where $B_{2n}$ denotes Bernoulli numbers. When $N$ is large, by taking the logarithm, substituting $B_2 = \frac{1}{6}$, and neglecting the higher-order terms, we obtain
\begin{align}
&\ \log{N!}-N\log{N}+N \\
&= N\log{\frac{N+1}{N}}+\frac{1}{2}\log{(N+1)}-1+\frac{1}{2}\log{(2\pi)}+\frac{1}{12(N+1)}. \notag
\end{align}
Using (24), we derive the difference between the logarithm of $P_H$ from (20) and $N\cdot MI$ from (18) as
\begin{align}
 & ER_1  \\
                   &= \frac{1}{2}\log{(N+1)(NX_0+1)(NX_1+1)(NX_2+1)(NX_3+1)} \notag \\
                       &\ -\frac{1}{2}\log{[N(X_0+X_1)+1][N(X_0+X_2)+1][N(X_1+X_3)+1][N(X_2+X_3)+1]}\notag \\
                       &\ +N\log{\frac{N+1}{N}}\notag \\
                       &\ +\sum_{k=0}^{3}NX_k\log{\frac{NX_k+1}{NX_k}}\notag \\
                       &\ -N(X_0+X_1)\log{\frac{N(X_0+X_1)+1}{N(X_0+X_1)}}-N(X_0+X_2)\log{\frac{N(X_0+X_2)+1}{N(X_0+X_2)}}\notag \\
                       &\ -N(X_1+X_3)\log{\frac{N(X_1+X_3)+1}{N(X_1+X_3)}}-N(X_2+X_3)\log{\frac{N(X_2+X_3)+1}{N(X_2+X_3)}}\notag \\
                  &\ -1+\frac{1}{2}\log{(2\pi)}+\frac{1}{12(N+1)}+\sum_{k=0}^{3}\frac{1}{12(NX_k+1)}-\frac{1}{12(NX_0+NX_1+1)}\notag \\
                       &\ -\frac{1}{12(NX_0+NX_2+1)}-\frac{1}{12(NX_1+NX_3+1)}-\frac{1}{12(NX_2+NX_3+1)}. \notag
\end{align}

We first evaluate each term on the right-hand side of (25) from above. The sum of the first and second lines on the right-hand side is less than $\frac{1}{2}\log{(N+1)}$. The third line is
\begin{equation}
N\log{\frac{N+1}{N}} = N\log{\Big(1+\frac{1}{N}\Big)}<N\times \frac{1}{N}=1.
\end{equation}
We evaluate the fourth to sixth lines as
\begin{align}
&\ NX_0\log{\frac{NX_0+1}{NX_0}\times\frac{N(X_0+X_1)}{N(X_0+X_1)+1}\times\frac{N(X_0+X_2)}{N(X_0+X_2)+1}} \\
&< NX_0\log{\Bigg(1+\frac{1}{NX_0}\Bigg)}<1, \notag
\end{align}
and so on. Moreover,
\begin{equation}
\frac{1}{12(NX_0+1)}+\sum_{k=1}^{3}\frac{1}{12(NX_k+1)} <\frac{1}{24}+3\times \frac{1}{12}=\frac{7}{24}.
\end{equation}
Thus, we evaluate $ER_1$ from above as
\begin{align}
 ER_1 &< \frac{1}{2}\log{(N+1)}+5-1+\frac{1}{2}\log{(2\pi)}+\frac{1}{24}+\frac{7}{24} \\
                      &< \frac{1}{2}\log{(N+1)}+5.253. \notag
\end{align}

Next, we evaluate the right-hand side of (25) from below. Using
\begin{equation}
N\log{(1+\frac{1}{N})}>N(\frac{1}{N}-\frac{1}{2N^2})=1-\frac{1}{2N}
\end{equation}
and similar inequalities, we obtain
\begin{align}
& ER_1 \\
&> \frac{1}{2}\log{\frac{2(N+1)}{(N+1)^4}}+1-\frac{1}{2N}+\sum_{k=0}^{3}(1-\frac{1}{2NX_k})-5+\frac{1}{2}\log{(2\pi)}+\frac{1}{12(N+1)}\notag \\
&> -\frac{3}{2}\log{(N+1)}+\frac{1}{2}\log{2}-\frac{1}{2N}-2+\frac{1}{2}\log{(2\pi)}+\frac{1}{12(N+1)}\notag \\
&> -\frac{3}{2}\log{(N+1)}+\frac{1}{12(N+1)}-\frac{1}{2N}-0.735.\notag
\end{align}
Thus, from (29) and (31),
\begin{align}
-\frac{3}{2}\log{(N+1)}+\frac{1}{12(N+1)}-\frac{1}{2N}-0.735 
&< ER_1 \\
&< \frac{1}{2}\log{(N+1)}+5.253. \notag
\end{align}
Hence, we have evaluated the first part of the error $ER_1$ from both above and below.

\subsubsection{Evaluation of the second part of the error}
Next, we estimate the second part of the error
\[
ER_2=-\log P_H-(-\log P_F)=\log\frac{P_F}{P_H},
\]
which quantifies the contribution of the Fisher tail sum relative to the single-table probability $P_H$.

Fix the margins in Table 1. The set of feasible $2\times2$ tables under these margins is finite.
Equivalently, letting $T$ denote the count in the $(A_1,B_1)$ cell, $T$ ranges over an integer interval $[T_{\min},T_{\max}]$.
Hence the number of feasible tables is
\[
M_N:=T_{\max}-T_{\min}+1\le N+1.
\]

By definition of Fisher's exact-test $p$-value in (12), we sum probabilities of feasible tables whose hypergeometric probability does not exceed that of the observed table.
Therefore each summand is at most $P_H$, and there are at most $M_N$ such tables.
Hence,
\begin{equation}
P_H \le P_F \le M_N P_H \le (N+1)P_H.
\end{equation}
Consequently,
\begin{equation}
0\le ER_2 \le \log(N+1).
\end{equation}
This two-sided bound is crude but sufficient for logarithmic asymptotics, since it contributes only an $O(\log N)$ term to $-\log P_F$.

\subsubsection{Combining the two parts of the error}

By combining the above bounds on $ER_1$ and $ER_2$, we can evaluate the total error as
\vspace{3pt}

\begin{align}
&\ -\log(N+1)-\frac{3}{2}\log{(N+1)}+\frac{1}{12(N+1)}-\frac{1}{2N}-0.735 \\
&< -\log{P_F}-N\cdot MI = ER_1-ER_2 \notag \\
&< \frac{1}{2}\log{(N+1)}+5.253. \notag
\end{align}
Therefore, transposing $N\cdot MI$,
\vspace{3pt}

\begin{align}
&\ N\cdot MI-\log(N+1)-\frac{3}{2}\log{(N+1)}+\frac{1}{12(N+1)}-\frac{1}{2N}-0.735 \\
&< -\log{P_F}\notag \\
           &< N\cdot MI +\frac{1}{2}\log{(N+1)}+5.253. \notag
\end{align}
Dividing by $N$, we obtain
\begin{align}
& MI+\frac{1}{N}\Bigg[-\log(N+1)-\frac{3}{2}\log{(N+1)}+\frac{1}{12(N+1)}-\frac{1}{2N}-0.735\Bigg] \\
&< -\frac{1}{N}\log{P_F}\notag \\
&< MI+\frac{1}{2N}\log{(N+1)}+\frac{5.253}{N}.\notag
\end{align}
\vspace{5pt}
As $N$ goes to infinity, we obtain
\begin{equation}
-\frac{1}{N}\log P_F = MI + O\!\left(\frac{\log(N+1)}{N}\right),
\end{equation}
\vspace{5pt}
because the terms on the left-hand and right-hand sides, except $MI$, converge to $0$. \\
Hence, we have proved Theorem 3 for a $2 \times 2$ contingency table. $\square$

\vspace{20pt}

\section{Proof of Theorem 3 for a general $m \times n$ contingency table}

In this section, we prove Theorem \ref{th3} for a general $m \times n$ contingency table. We consider the following two $m \times n$ contingency tables, Tables 3 and 4, for two random variables $A$ and $B$. $A$ takes values from $A_1$ to $A_m$ and $B$ takes values from $B_1$ to $B_n$. Table 3 shows the frequency of the combination of variables. Meanwhile, Table 4 shows the relative frequency obtained from Table 3 by dividing by $N$. Similar to Tables 1 and 2, Tables 3 and 4 can be converted to each other. Let $X_{ij}$ be the relative frequency by which $A$ takes $A_i$ and $B$ takes $B_j$.

\vspace{10pt}

\begin{center}
T{\small ABLE} 3\\
\textit{$m\times n$ contingency table of the observed frequency.} \\
\textit{All the values of the cells are non-negative integers.} \\
\begin{tabular}{ccccccc} \hline
\rule[-0.5em]{0em}{2em} & $A_1$ & $\dots$ & $A_i$ & $\dots$ & $A_m$ & Total\\ [3pt]\hline
\rule[-0.5em]{0em}{2em} $B_1$ & $NX_{11}$ & $\dots$ & $NX_{i1}$ & $\dots$ & $NX_{m1}$ & $N\sum_{i=1}^{m} X_{i1}$ \\  [3pt]
$\dots$ & $\dots$ & $\dots$ & $\dots$ & $\dots$ & $\dots$ & $\dots$ \\ [3pt]
$B_j$ & $NX_{1j}$ & $\dots$ & $NX_{ij}$ & $\dots$ & $NX_{mj}$ & $N\sum_{i=1}^{m} X_{ij}$ \\ [3pt]
$\dots$ & $\dots$ & $\dots$ & $\dots$ & $\dots$ & $\dots$ & $\dots$ \\ [3pt]
$B_n$ & $NX_{1n}$ & $\dots$ & $NX_{in}$ & $\dots$ & $NX_{mn}$ & $N\sum_{i=1}^{m} X_{in}$ \\ [3pt]\hline
\rule[-0.5em]{0em}{2em} Total & $N\sum_{j=1}^{n} X_{1j}$ & $\cdots$ & $N\sum_{j=1}^{n} X_{ij}$ & $\cdots$ & $N\sum_{j=1}^{n} X_{mj}$ & $N$ \\ [3pt]\hline
\end{tabular}
\end{center}

\vspace{10pt}

\begin{center}
T{\small ABLE} 4\\
\textit{$m\times n$ contingency table of the relative frequency.}  \\
\textit{All the values of the cells are non-negative real numbers.}
\begin{tabular}{cccccccc} \hline
\rule[-0.5em]{0em}{2em} &  $A_1$ & $\dots$ & $A_i$ & $\dots$ & $A_m$ & & Total\\ [3pt]\hline
\rule[-0.5em]{0em}{2em} $B_1$ & $X_{11}$ & $\dots$ & $X_{i1}$ & $\dots$ & $X_{m1}$ & & $\sum_{i=1}^{m} X_{i1}$ \\ [3pt]
$\dots$ & $\dots$ & $\dots$ & $\dots$ & $\dots$ & $\dots$ & & $\dots$ \\ [3pt]
$B_j$ & $X_{1j}$ & $\dots$ & $X_{ij}$ & $\dots$ & $X_{mj}$ & & $\sum_{i=1}^{m} X_{ij}$ \\ [3pt]
$\dots$ & $\dots$ & $\dots$ & $\dots$ & $\dots$ & $\dots$ & & $\dots$  \\ [3pt]
$B_n$ & $X_{1n}$ & $\dots$ & $X_{in}$ & $\dots$ & $X_{mn}$ & & $\sum_{i=1}^{m} X_{in}$ \\ [3pt]\hline
\rule[-0.5em]{0em}{2em} Total &  $\sum_{j=1}^{n} X_{1j}$ & $\dots$ & $\sum_{j=1}^{n} X_{ij}$ & $\dots$ & $\sum_{j=1}^{n} X_{mj}$ & & $1$ \\ [3pt]\hline
\end{tabular}
\end{center}
\vspace{15pt}

Using Tables 3 and 4, we calculate and evaluate $-\log {P_H} - N\cdot MI$ from above similar to (25) to (29) as

\begin{align}
&\ -\log{P_H}-N\cdot \MI \\
&= \frac{1}{2}\log {\Bigg[(N+1)\prod_{i=1}^m\prod_{j=1}^n(NX_{ij}+1)\Bigg]}\notag \\
&\ -\frac{1}{2}\log {\Bigg[\prod_{j=1}^n(N\sum_{i=1}^mX_{ij}+1)\Bigg]\Bigg[\prod_{i=1}^m(N\sum_{j=1}^nX_{ij}+1)\Bigg]}\notag \\
&\ +N\log{\frac{N+1}{N}}+\sum_{i=1}^m\sum_{j=1}^nNX_{ij}\log{\frac{NX_{ij}+1}{NX_{ij}}}\notag \\
&\ -N\sum_{j=1}^n\Bigg(\sum_{i=1}^mX_{ij}\Bigg)\log{\frac{N\sum_{i=1}^mX_{ij}+1}{N\sum_{i=1}^mX_{ij}}}\notag \\
&\ -N\sum_{i=1}^m\Bigg(\sum_{j=1}^nX_{ij}\Bigg)\log{\frac{N\sum_{j=1}^nX_{ij}+1}{N\sum_{j=1}^nX_{ij}}}\notag \\
&\ -1+\frac{1}{2}\log(2\pi)+\frac{1}{12(N+1)}+\sum_{i=1}^m\sum_{j=1}^n\frac{1}{12(NX_{ij}+1)}\notag \\
&\ -\sum_{j=1}^n\frac{1}{12\Bigg(N\sum_{i=1}^mX_{ij}+1\Bigg)}-\sum_{i=1}^m\frac{1}{12\Bigg(N\sum_{j=1}^nX_{ij}+1\Bigg)}\notag \\
&< \frac{1}{2}\log{(N+1)}+1+mn-1+\frac{1}{2}\log{(2\pi)}+\frac{1}{24}+\frac{1}{24}+\frac{1}{12}(mn-1)\notag \\
&< \frac{1}{2}\log{(N+1)}+\frac{13}{12}mn+0.919. \notag
\end{align}
\vspace{5pt}

Meanwhile, we evaluate $-\log {P_H} - N\cdot MI$ from below similar to (30) and (31) as
\begin{align}
&\ -\log{P_H}-N\cdot MI \\
&> \frac{1}{2}\log{\frac{2(N+1)}{(N+1)^{mn}}}+1-\frac{1}{2N}+\sum_{i=1}^m\sum_{j=1}^n\Bigg(1-\frac{1}{2NX_{ij}}\Bigg)\notag \\
&\ -mn-1+\frac{1}{2}\log{(2\pi)}+\frac{1}{12(N+1)}\notag \\
&> -\frac{mn-1}{2}\log{(N+1)}+\frac{1}{2}\log{2}-\frac{1}{2N}-\frac{mn}{2}+\frac{1}{2}\log{(2\pi)}+\frac{1}{12(N+1)}\notag \\
&> -\frac{mn-1}{2}\log{(N+1)}+\frac{1}{12(N+1)}-\frac{1}{2N}-\frac{mn}{2}+1.265. \notag
\end{align}
\vspace{5pt}

Finally, we evaluate the difference between $P_F$ and $P_H$ in the two-sided definition (12).
Let $\mathcal{F}_N$ denote the finite set of feasible $m\times n$ contingency tables with the fixed margins.
Since every summand in (12) is at most $P_H$, and the number of feasible tables is $|\mathcal{F}_N|$, we have
\begin{equation}
P_H \le P_F \le |\mathcal{F}_N|\,P_H.
\end{equation}
A crude but convenient bound is $|\mathcal{F}_N|\le (N+1)^{(m-1)(n-1)}$, since an $m\times n$ table with fixed margins is determined by $(m-1)(n-1)$ free entries and each entry lies in $\{0,1,\ldots,N\}$.
Therefore,
\begin{equation}
P_H \le P_F \le (N+1)^{(m-1)(n-1)}P_H.
\end{equation}
This bound is crude but sufficient for establishing logarithmic asymptotics, since it contributes only a polynomial factor on the probability scale and hence vanishes on the $(1/N)\log$-scale.

Hence,
\begin{equation}
-\log P_H-(m-1)(n-1)\log(N+1)-N\cdot MI < -\log P_F-N\cdot MI \le -\log P_H-N\cdot MI.
\end{equation}

Combining the lower bound on $-\log P_H-N\cdot MI$ derived above with this inequality, we obtain
\begin{align}
&\ -\frac{mn-1}{2}\log{(N+1)}+\frac{1}{12(N+1)}-\frac{1}{2N}-\frac{mn}{2}+1.265 \\
& -(m-1)(n-1)\log(N+1)\notag \\
&< -\log{P_F}-N\cdot MI. \notag
\end{align}
Meanwhile, using the upper bound on $-\log P_H-N\cdot MI$, we have
\begin{equation}
-\log{P_F}-N\cdot MI < -\log{P_H}-N\cdot MI < \frac{1}{2}\log{(N+1)}+\frac{13}{12}mn+0.919.
\end{equation}
Therefore,
\begin{align}
&\ N\cdot MI -\frac{mn-1}{2}\log{(N+1)}+\frac{1}{12(N+1)}-\frac{1}{2N} \\
&\ -\frac{mn}{2}+1.265-(m-1)(n-1)\log(N+1)\notag \\
&< -\log{P_F}  \notag \\
&< N\cdot MI+\frac{1}{2}\log{(N+1)}+\frac{13}{12}mn+0.919. \notag
\end{align}

Hence, by dividing by $N$,
\begin{equation}
-\frac{1}{N}\log P_F = MI + O\!\left(\frac{\log(N+1)}{N}\right),
\end{equation}
as $N$ goes to infinity. Thus, we have proved Theorem \ref{th3} for an $m \times n$ contingency table. $\square$

\vspace{10pt}

\section{Numerical simulations}
\begin{figure}
\centering
\includegraphics[width=\columnwidth]{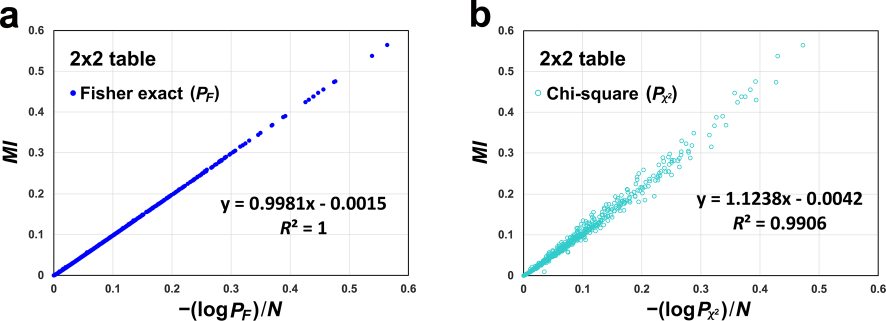}
\caption{\textit{$P_F$, $P_{\chi^2}$ and $MI$ of $2\times 2$ contingency tables. (a) $P_F$ and $MI$. (b) $P_{\chi^2}$ and $MI$. The equations and $R^2$ represent the regression lines and the determination coefficients between $-(\log{P_F})/N$ and $MI$, and between $-(\log{P_{\chi^2}})/N$ and $MI$, respectively.}}
\label{Figure 2}
\end{figure}

To assess the validity of Theorem \ref{th3}, we performed Monte Carlo simulations. We created $2\times 2$ contingency tables using random numbers. Because we are interested in the sample size for which Fisher's exact test is cumbersome, we set the sample size $N$ of the tables to 1000, which is practical in statistical analysis in medicine and biology. We ran this trial 1000 times, and then calculated $P_F$, $MI$ and the chi-square test $p$-value ($P_{\chi^2}$) of the tables. It is well-known that $P_{\chi^2}$ provides a good approximation to $P_F$ when $N$ is sufficiently large. We wrote the computer programs in MATLAB. As shown in Figure \ref{Figure 2}(a), when we took the logarithm of $P_F$ to observe the relationship between $P_F$ and $MI$, $-\log {P_F}$ and $MI$ were scattered exactly along the line $MI = -(\log {P_F})/N$ ($R^2=1$). Meanwhile, when we took the logarithm of $P_{\chi^2}$ to observe the relationship between $P_{\chi^2}$ and $MI$, $-\log{P_{\chi^2}}$ and $MI$ were more scattered around the line ($R^2=0.9906$), as shown in Figure \ref{Figure 2}(b).

\begin{figure}
\centering
\includegraphics[width=\columnwidth]{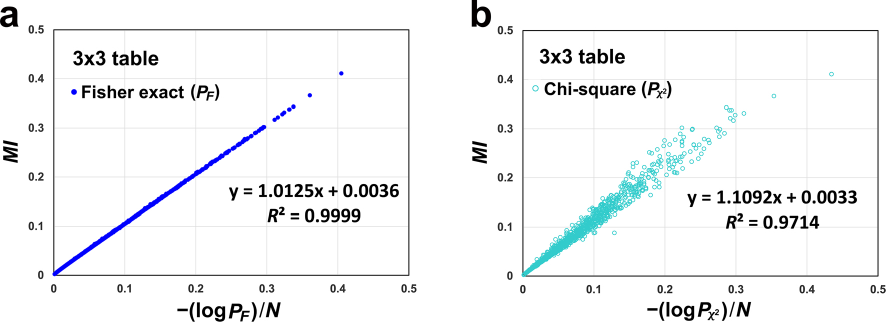}
\caption{\textit{$P_F$, $P_{\chi^2}$ and $MI$ of $3\times 3$ contingency tables. (a) $P_F$ and $MI$. (b) $P_{\chi^2}$ and $MI$. The equations and $R^2$ represent the regression lines and the determination coefficients between $-(\log{P_F})/N$ and $MI$, and between $-(\log{P_{\chi^2}})/N$ and $MI$, respectively.}}
\label{Figure 3}
\end{figure}

Next, we conducted similar experiments by creating $3\times 3$ contingency tables. In this case, we calculated $P_F$ using the statistical software Stata. Then, $-\log{P_F}$ and $MI$ followed the equation $MI=-(\log{P_F})/N$ ($R^2=0.9999$), as indicated in Figure \ref{Figure 3}(a). Meanwhile, the $R^2$ of $-\log{P_{\chi^2}}$ and $MI$ was 0.9714, as shown in Figure \ref{Figure 3}(b). Therefore, $P_F$ and $MI$ were converted to each other and Theorem \ref{th3} was true for both in $2\times 2$ and $3\times 3$ contingency tables. By contrast, the correlation of $P_{\chi^2}$ and $MI$ was worse, with a poorer coefficient of determination than $P_F$ and $MI$.

\section{Meta-analysis}
Meta-analysis integrates data from multiple lines of studies \cite{r3,r8}, thereby yielding reliable statistics, often with a decreased $p$-value. It has been applied to estimate MI \cite{r5, r15, r18}, where the authors computed the weighted average of MI without producing its $p$-value. By contrast, our method enables us to calculate the $p$-value of MI, which allows the application of meta-analysis to this research direction from a novel point of view.
\vspace{2pt}

\subsection{Integration of contingency tables with the same measurement error}
First, we mention the case in which the data regard the same random variables and integrate $H$ contingency tables represented by Table 5 similarly to Table 1, where $H$ is sufficiently large. This approach applies to the case in which the measurement errors of the $H$ contingency tables are the same. In Table 5, the index $h$ runs from 1 to $H$, $N_h$ is the sample size of the $h$-th table, and $X_{hi}$ ($i$ = 0, 1, 2, 3) are the observed relative frequencies in the $h$-th table.

We integrate the $H$ tables represented by Table 5 into Table 6 by summing each cell. In Table 6, we perform the summation with respect to $h$ from 1 to $H$. Then we divide the cells of Table 6 by $N_s=\sum_{h=1}^H N_h$ to obtain Table 7, in which $Z_i = \sum_{h=1}^H N_hX_{hi} / N_s$ ($i$ = 0, 1, 2, 3). If we calculate $P_F$ from Table 6 and $MI$ from Table 7, we obtain asymptotically
\begin{equation}
MI_s=-\frac {1}{N_s} \log{P_F},
\end{equation}
where $MI_s$ is the integrated $MI$ calculated from Table 7.

\vspace{3pt}
This formula can be proved similarly to the proof of Theorem \ref{th3} for $2\times 2$ contingency tables. The same formula holds for general $m\times n$ contingency tables. Thus, Theorem \ref{th3} can also be applied to meta-analysis to enable us to perform a more accurate estimation of $MI$. (48) demonstrates that, although the observed MI may differ table by table depending on $h$, we can estimate the true MI precisely using Tables 6 and 7 of the meta-analysis according to the law of large numbers applied to $Z_0$ to $Z_3$. Additionally, we expect the $P_F$ of Table 6 to be less than or equal to those of Table 5 because the number of data $N_s$ should be sufficiently large, which guarantees a more accurate estimation of $MI$.
\newpage

\begin{center}
T{\small ABLE} 5\\
\textit{Contingency table of the observed frequency similar to Table 1.} \\
\textit{All the values of the cells are non-negative integers.} \\
\begin{tabular}{cccc} \hline
\rule[-0.5em]{0em}{2em} & $A_1$ & $A_2$ & Total\\ [3pt]\hline
\rule[-0.5em]{0em}{2em} $B_1$ & $N_hX_{h0}$ & $N_hX_{h1}$ & $N_h(X_{h0}+X_{h1})$ \\  [3pt]
$B_2$ & $N_hX_{h2}$ & $N_hX_{h3}$ & $N_h(X_{h2}+X_{h3})$ \\ [3pt]\hline
\rule[-0.5em]{0em}{2em} Total & $N_h(X_{h0}+X_{h2})$ & $N_h(X_{h1}+X_{h3})$ & $N_h$ \\ [3pt]\hline
\end{tabular}
\end{center}

\vspace{20pt}

\begin{center}
T{\small ABLE} 6\\
\textit{Contingency table of the observed frequency integrating the $H$ tables represented in Table 5.} \\
\textit{All the values of the cells are non-negative integers.} \\
\begin{tabular}{cccc} \hline
\rule[-0.5em]{0em}{2em} & $A_1$ & $A_2$ & Total\\ [3pt]\hline
\rule[-0.5em]{0em}{2em} $B_1$ & $\sum N_hX_{h0}$ & $\sum N_hX_{h1}$ & $\sum N_h(X_{h0}+X_{h1})$ \\  [3pt]
$B_2$ & $\sum N_hX_{h2}$ & $\sum N_hX_{h3}$ & $\sum N_h(X_{h2}+X_{h3})$ \\ [3pt]\hline
\rule[-0.5em]{0em}{2em} Total & $\sum N_h(X_{h0}+X_{h2})$ & $\sum N_h(X_{h1}+X_{h3})$ & $\sum N_h$ \\ [3pt]\hline
\end{tabular}
\end{center}
\vspace{20pt}

\begin{center}
T{\small ABLE} 7\\
\textit{Contingency table of the relative frequency obtained from Table 6.}  \\
\textit{All the values of the cells are non-negative real numbers.}
\begin{tabularx}{90mm}{cccccccc} \hline
\rule[-0.5em]{0em}{2em} & & $A_1$ & & $A_2$ & & & Total\\ [3pt]\hline
\rule[-0.5em]{0em}{2em} $B_1$ & & $Z_0$ & & $Z_1$ & & & $Z_0+Z_1$ \\ [3pt]
$B_2$ & & $Z_2$ & & $Z_3$ & & &  $Z_2+Z_3$ \\ [3pt]\hline
\rule[-0.5em]{0em}{2em} Total & & $Z_0+Z_2$ & & $Z_1+Z_3$ & & & $1$ \\ [3pt]\hline
\end{tabularx}
\end{center}
\vspace{10pt}

\subsection{Integration of contingency tables with different measurement errors}
Finally, we extend (48) to the case in which the contingency tables are represented by the same form as Table 5 but observe features of the same random variables with different measurement errors, and we can again refine the $p$-value of MI. In this case, there is no prior information about the distributions of the random variables. Then it is appropriate to apply the maximum entropy principle \cite{r9}.

Let $MI_h$ be the observed MI in the $h$-th table represented by Table 5 and let $MI_s$ be MI obtained by integrating multiple tables. Because both margins are not known exactly in advance, it is desirable to use measures other than the $P_F$ of the tables at first \cite{r1}. Let $p_h$ be the $p$-value estimated from $MI_h$. Then, from (8),
\begin{equation}
N_hMI_h = -\log{p_h}.
\end{equation}
If the $H$ tables represented by Table 5 are the results of independent trials, then
\begin{equation}
MI_s=\frac{1}{N_s}\sum_{h=1}^H N_hMI_h.
\end{equation}
$N_hMI_h$ can be evaluated using $-\log{P_{Fh}}$ from above and below according to (37), where $P_{Fh}$ is the $P_F$ of the $h$-th table. Therefore, when $N_s$ is sufficiently large, $\sum_{h=1}^H N_hMI_h/N_s$ can be replaced by $-(\log{\prod_{h=1}^H P_{Fh}})/N_s$.
Then, asymptotically
\begin{equation}
MI_s = -\frac{1}{N_s}\log{\prod_{h=1}^H P_{Fh}} = -\frac{1}{N_s}\log{p_{s}},
\end{equation}
where $p_s$ is the $p$-value of $MI_s$ under the maximum entropy principle. (50) means that $MI_s$ is the weighted average of $MI_h$, and approaches the true MI when $N_s$ becomes sufficiently large. Additionally, $p_s$ is less than or equal to every $P_{Fh}$ when $H$ is large, which demonstrates that the $p$-value of MI is again refined by meta-analysis. Moreover, the same can be said of general $m\times n$ contingency tables.

In our methods in 6.2, the weight for each $MI_h$ is the sample size, $N_h$. By contrast, in \cite{r5, r15, r18}, the weighted average of MI was computed, where the weights used were the inverse of the $MI$ variance, which is related to sample size. Despite this, those methods could not calculate the $p$-value of MI. The advantage of our method is that it allows a facile calculation of the $p$-value of the weighted average of MI. Thus, from the two complementary conditions of 6.1 and 6.2, we have demonstrated that our theorems can refine MI using meta-analysis, thereby producing more decreased $p$-values.

\section{Discussion and application}
This paper establishes logarithmic (exponential-rate) relations between statistical significance in testing independence and information-theoretic dependence quantified by Shannon mutual information (MI), in two complementary settings.
When no prior distributional information is available, the maximum-entropy principle yields an exponential-form calibration for the probability of realizing a given magnitude of information exchange (Theorems~1--2).
When marginal totals are fixed, we analyze Fisher's exact test and show that its two-sided $p$-value $P_F$ satisfies
\[
-\frac{1}{N}\log P_F = MI + O\!\left(\frac{\log(N+1)}{N}\right),
\]
so that, equivalently,
\[
P_F=\exp\{-N\cdot MI+O(\log N)\},
\]
with explicit finite-$N$ bounds derived in the proofs.
Taken together, these results clarify how MI governs the exponential rate at which $p$-values decay with sample size in the settings studied here.

\vspace{6pt}
\noindent
We highlight several implications and practical uses.

\vspace{5pt}
\noindent
\textbf{(i) Computational surrogate for large tables.}
Exact computation of Fisher's $p$-value can be demanding for large sample sizes or for larger $m$ and $n$ (even when $m$ and $n$ are moderate), because the $p$-value is a tail sum over many feasible tables \cite{r12}.
Theorem~\ref{th3} shows that, on the logarithmic scale and for fixed $m,n$, $-(1/N)\log P_F$ is well-approximated by $MI$, with a discrepancy of order $\log N/N$.
Thus, $MI$ provides a fast proxy for the \emph{exponential rate} of $P_F$ when $N$ is large, which may be useful for screening or for benchmarking dependence across many tables.

\vspace{5pt}
\noindent
\textbf{(ii) Information-theoretic interpretation of significance.}
The relations in Theorems~1--3 imply that $-(1/N)\log(p\text{-value})$ can be interpreted as an information-like quantity (in nats per observation).
In particular, smaller $p$-values correspond to larger information exchange, while the factor $N$ explains how evidence accumulates with sample size.

\vspace{5pt}
\noindent
\textbf{(iii) Comparing dependence across sample sizes.}
A well-known feature of classical $p$-values is that, for any fixed nonzero dependence, $p$-values tend to decrease as $N$ increases.
Our results make this scaling explicit: in the fixed-margin setting, $P_F$ decays roughly like $\exp(-N\cdot MI)$ up to a polynomial factor.
This reinforces the role of $MI$ as a stable effect-size measure for comparing dependence across studies with different sample sizes, while still allowing significance to be recovered on the exponential scale.

\vspace{5pt}
\noindent
\textbf{(iv) Model-based $p$-values for MI under maximum entropy.}
Mutual information is often used as a dependence measure but is not itself a significance statement.
In the maximum-entropy baseline setting, Theorem~\ref{th2} yields the explicit calibration $P_{N\cdot MI}=e^{-N\cdot MI}$, providing a direct model-based tail-probability interpretation of an observed $MI$.

\vspace{5pt}
\noindent
\textbf{(v) Combining evidence via meta-analysis.}
Because the logarithm of a $p$-value is additive across independent studies, the exponential-rate viewpoint suggests natural ways to combine evidence from multiple datasets.
In Section~6 we illustrated two such strategies: pooling contingency tables with comparable measurement conditions (Section~6.1) and combining studies with heterogeneous measurement error via a meta-analytic conversion on the probability scale (Section~6.2).
The present results justify these procedures on the logarithmic scale by linking $-\log(p\text{-value})$ to $N\cdot MI$ up to $O(\log N)$ corrections.

\vspace{8pt}
\noindent
\textbf{Scope and limitations.}
Theorem~\ref{th3} is proved for fixed $m,n$ as $N\to\infty$; understanding regimes where $m$ and/or $n$ grow with $N$, or obtaining sharp polynomial prefactors beyond the logarithmic scale, are natural directions for future work.

\begin{acks}[Acknowledgments]
We are grateful to T. Okada for valuable comments. We thank A. Goto for helpful advice about numerical examples. We thank Edanz \\ 
(https://www.edanz.com/ac) for editing a draft of this manuscript.
\end{acks}

\end{document}